\documentclass[]{interact}
 \usepackage[colorlinks=true]{hyperref}
\usepackage{epstopdf}
\usepackage{subfigure}

\usepackage[numbers,sort&compress]{natbib}
\bibpunct[, ]{[}{]}{,}{n}{,}{,}
\makeatletter
\def\NAT@def@citea{\def\@citea{\NAT@separator}}
\makeatother

\theoremstyle{plain}
\newtheorem{theorem}{Theorem}[section]
\newtheorem{lemma}[theorem]{Lemma}

\theoremstyle{definition}

\theoremstyle{remark}
\newtheorem{remark}{Remark}

\newtheorem{algorithm}[theorem]{Algorithm}
\begin{document}


\title{Golden ratio algorithms with new stepsize rules for variational inequalities}

\author{
\name{Dang Van Hieu\textsuperscript{a}, Yeol Je Cho\textsuperscript{b,c} and Yi-bin Xiao\textsuperscript{c}\thanks{Dang Van Hieu. Email: dangvanhieu@tdtu.edu.vn}}
\affil{\textsuperscript{a}Applied Analysis Research Group, Faculty of Mathematics and Statistics, Ton Duc Thang University, Ho Chi Minh City, Vietnam; \\
\textsuperscript{b}Department of Mathematics Education, Gyeongsang National University, Jinju 52828, Korea; \\
\textsuperscript{c}School of Mathematical Sciences, University of Electronic Science and Technology of China, Chengdu, Sichuan, 611731, P.R. China\\[.1in]
\textit{Dedicated to Professor Pham Ky Anh on the Occasion of his 70th Birthday}
}
}

\maketitle

\begin{abstract}
In this paper, we introduce two golden ratio algorithms with new stepsize rules for solving pseudomonotone and Lipschitz variational inequalities in 
finite dimensional Hilbert spaces. The presented stepsize rules allow the resulting algorithms to work without the prior knowledge of the Lipschitz 
constant of operator. The first algorithm uses a sequence of stepsizes which is previously chosen, diminishing and non-summable. While the stepsizes 
in the second one are updated at each iteration and by a simple computation. A special point is that the sequence of stepsizes generated by 
the second algorithm is separated from zero. The convergence as well as the convergence rate of the proposed algorithms are established under 
some standard conditions. Also, we give several numerical results to show the behavior of the algorithms in comparisons with other algorithms.
\end{abstract}


\section{Introduction}\label{intro}

In this paper, we focus on the following {\it variational inequality problem} (shortly, {\bf (VIP)}):
$$
\mbox{Find}~x^*\in \Re^m~\mbox{such that}~\left\langle Fx^*,x-x^*\right\rangle+g(x)-g(x^*)\ge 0,\,\,\,\forall x\in \Re^m,
$$
where $g:\Re^m \to (-\infty,+\infty]$ is a proper convex lower semi-continuous function with the domain ${\rm dom}\,g=\left\{x\in \Re^m: g(x)<+\infty\right\}$ 
and $F: {\rm dom}\,g\to \Re^m$ is an operator. The function $g$ here cannot be smooth. Recall that the {\it proximal operator} ${\rm prox}_g$ of $g$ is defined by 
$$
{\rm prog}_g(x)=\arg\min_{y\in \Re^m} \Big\{g(y)+\frac{1}{2}||y-x||^2\Big\}.
$$
We are also interested here in the problem {\bf (VIP)} where the proximal mapping of $g$ is computable. The problem {\bf (VIP)} is known as a central problem 
in nonlinear analysis, especially, in optimization, control theory, games theory \cite{DGM2003,FP2002,GMP2004,BT2009,KS1980,M2018,K2000,K2007} and 
other fields \cite{ADJS18,SNAT18,A17,MA19,AJS18,RCAA17,SXC19,WXWC16,PPXY18}. Considering the problem {\bf (VIP)} in a special case, when 
$g=\delta_C$, the indicator operator of a nonempty closed convex set $C$ in $\Re^m$, this problem reduces to the classical {\it variational inequality problem} 
\cite{KS1980,S1964}:
$$
\mbox{Find}~x^*\in C~\mbox{such that}~\left\langle Fx^*,x-x^*\right\rangle\ge 0,\,\,\,\forall x\in C.
\eqno{\rm \bf (VIP)}
$$
Moreover, the motivations of studying the problem {\bf (VIP)} come from optimization point of view. Several models arising naturally can be formulated 
as the problem {\bf (VIP)}, see, e.g, in \cite{FP2002,KS1980,M2018}. 
We restrict our interest in the following two problems. The first basic problem is a {\it convex-concave saddle point problem}:
\begin{equation}\label{ccspp}
\min_{x\in \Re^m} \max_{y\in \Re^n} L(x,y) = g_1(x) +K(x,y) - g_2(y),
\end{equation}
where $x\in \Re^m$, $y\in \Re^n$, $g_1:\Re^m\to (-\infty,+\infty]$, $g_2:\Re^n\to (-\infty,+\infty]$ are proper convex lower semi-continuous functions, 
$K:{\rm dom}\,g_1 \times {\rm dom}\,g_2\to \Re$ is a smooth convex-concave function. The saddle point problem (\ref{ccspp}) can be considered the problem {\bf (VIP)} with 
$g(z)=g_1(x)+g_2(y)$, $F(z)=\left[\nabla_x K(x,y); -\nabla_y K(x,y)\right]^T$ and $z=(x,y)\in \Re^{m+n}$. This is a type of example for nonsmooth problem, 
where the gradient method (one step) \cite{AH1958} cannot work. The early proposed iterative methods for this problem may be the extragradient methods (two steps) 
\cite{K1976,P1980}. Recently, many works on the problem {\bf (VIP)} and related problems  have been devoted to proposing different projection-like methods under various types of conditions, 
such as the subgradient extragradient methods \cite{CGR2011JOTA,CGR2011OMS,CGR2012,HT2018,H2017COAP}, the modified extragradient method 
\cite{MS2014}, the projected reflected gradient method \cite{MG2016,M2015} and others \cite{C2008,D2018,HAM2017,KRS2011,KV2014,M2008,T2000,XAY2017,Y2009,Y2010}. 
\vskip 1mm

Another model arising in the signal processing literature is a nonsmooth convex optimization model,
\begin{equation}\label{ncom}
\min_{x\in \Re^m} J(x)=f(x)+g(x),
\end{equation}
where $x\in \Re^m$ and $f,~g:\Re^m\to \Re$ with $g$ possibly nonsmooth. This model is equivalent to the problem {\bf (VIP)} with $F=\nabla f$. 
Applying generic variational methods to solve special models as the optimization problem (\ref{ncom}) is not a good choice. Optimization methods 
in general have better theoretical convergence rates because they can exploit the characteristics of the potential operator $\nabla f$. However, this seems 
only true when $\nabla f$ is Lipschitz continuous. Without such a condition, the proximal gradient methods cannot hold anymore. Recently, in the case 
with non-Lipschitz $\nabla f$, some notable optimization methods can be found in, for example, in \cite{BBT2016,LFN2018,TKC2015}. Especially, the so-called 
NoLips algorithm developed in \cite{BBT2016} is an interesting and promising method which can use a fixed stepsize when $\nabla f$ is non-Lipschitz continuous. 
However, the obtained results are not generic because they depend strictly on problem instances and used Bregman distances.
\vskip 1mm

Methods for solving the problem {\bf (VIP)} without the Lipschitz continuity of $F$ often use a linesearch proceduce which runs in each iteration 
of the algorithm until a stopping criterion is satisfied. Linesearch methods are thus time-consuming because they require 
many computations of values 
of $F$ as well as projections onto feasible set. Besides, the estimates of complexity in linesearch methods become not so informative. This is clear because 
they only show how many outer iterations needed to be done to obtain the desired accuracy while the number of inner linesearch iterations cannot be 
mentioned. Moreover, in the case $F$ is even Lipschitz continuous, but in general the Lipschitz constant is often unknown and in nonlinear problems, it can be 
difficult to approximate. Recently, without any linesearch procedure, some interesting methods for solving the classical variational inequalities can be 
found, for instance, in \cite{K1989,M2018,T2004} where stepsizes are updated over iteration by some cheap computations.
\vskip 1mm

Especially, an interesting idea has been developed recently by Malitsky in \cite{M2018}. He has proposed a new algorithm, named the 
{\it explicit golden ratio algorithm} (shortly, EGRAAL), for solving problem {\bf (VIP)} when $F$ is locally Lipschitz continuous. 
The EGRAAL has a simple and elegant structure and only requires the computations of one proximal mapping value (with function $g$) 
and one value of $F$. His algorithm is explicit in the sense that it does not require any linesearch procedure. Its stepsizes are computed 
explicitly at each iteration from the previous iterates. The theoretical and numerical results in \cite{M2018} are promising and also suggest 
some directions for studying in the future (see, \cite[Sects. 5 and 6]{M2018}).
\vskip 2mm

In this paper, motivated by the results in \cite{K1989,M2018,T2004}, we propose two different golden ratio algorithms with new 
stepsize rules for solving pseudomonotone problem {\bf (VIP)} with a Lipschitz condition. The stepsize strategies here are simpler than those ones in \cite{M2018}. 
The variable stepsizes in the new algorithms are chosen previously or computed easily. Also, the resulting algorithms work without any information of Lipschitz 
constant of operator, i.e., the Lipschitz constant must not be the input parameters of the algorithms. More precisely, we first consider a golden 
ratio algorithm with a priorly taken sequence of stepsizes being diminishing and non-summable. Following to this strategy, the algorithm works well for problem 
{\bf (VIP)} with the strong pseudomonotonicity of $F$. For a weaker assumption of pseudomonotonicity of $F$, we present 
the second golden ratio algorithm incorporated with variable stepsizes updating step-by-step. 
The convergence as well as the convergence rate of the proposed algorithms are established. Finally, the theoretical results are confirmed by several our 
numerical experiments in comparisons with other known algorithms.
\vskip 2mm

The remainder of this paper is organized as follows: In Section \ref{main}, we introduce a golden ratio algorithm with a sequence of stepsizes priorly chosen. 
Section \ref{main1} deals with another golden ratio algorithm with a simpler stepsize rule. Finally, in Section \ref{example}, we study the numerical behaviour of 
the new algorithms on two test problems and compare them with others.

\section{A golden ratio algorithm with diminishing stepsizes}\label{main}

In this section, we present a golden ratio algorithm with a sequence of stepsizes being diminishing and non-summable. The use of this new stepsize rule 
allows the algorithm to work without previously knowing the Lipschitz constant of the operator. We set $\varphi=\frac{\sqrt{5}+1}{2}$, which is called 
the {\it golden ratio}. Moreover, we denote $VI(F,g)$ the solution set of the problem 
{\bf (VIP)} and it is assumed to be nonempty. The following is the algorithm in details:
\vskip 2mm

\begin{algorithm}[Golden Ratio Algorithm with diminishing stepsizes].\label{alg1}\\
\textbf{Initialization:} Choose $x_1,~\bar{x}_0\in \Re^m$ and a non-increasing sequence of stepsizes $\left\{\lambda_n\right\}\subset (0,+\infty)$ such that 
the following conditions hold:
$$
\mbox{\rm (H1)}\,\,\,\lim_{n\to \infty}\lambda_n=0,\qquad \mbox{\rm (H2)}\,\,\,\sum_{n=1}^\infty \lambda_n=+\infty,\qquad \mbox{\rm (H3)}\,\,\,\liminf_{n\to \infty}\frac{\lambda_n}{\lambda_{n-1}}>0.
$$
\textbf{Iterative Steps:} Assume that $x_n,~\bar{x}_{n-1}$ are known, calculate $x_{n+1}$ as follows:
$$
\left \{
\begin{array}{ll}
\bar{x}_n=\frac{(\varphi-1)x_n+\bar{x}_{n-1}}{\varphi},\\
x_{n+1}=\mbox{\rm prox}_{\lambda_n g}(\bar{x}_n-\lambda_n F(x_n)).
\end{array}
\right.
$$
\end{algorithm}
\vskip 2mm

\noindent An example for the sequence $\left\{\lambda_n\right\}$ satisfying the conditions (H1)-(H3) is $\lambda_n=\frac{1}{n^p}$ with $0<p\le 1$. In order to 
establish the convergence of Algorithm \ref{alg1}, we assume that the operator $F:{\rm dom}g\to \Re^m$ satisfies the following conditions:
\vskip 1mm

(SP)\, $F$ is strongly pseudomonotone, i.e., there exists $\gamma>0$ such that 
$$
 \left\langle F(x),x-y\right\rangle \ge 0 \Longrightarrow  \left\langle F(y),x-y\right\rangle \ge \gamma ||x-y||^2,\,\,\,\forall x,~y\in \Re^m;
$$

(LC)\, $F$ is Lipschitz continuous, i.e., there exists $L>0$ such that 
$$
 ||F(x)-F(y)|| \le L||x-y||,\,\,\,\forall x,~y\in \Re^m.
$$
\vskip 1mm

However, it is not necessary to know the two constants $\gamma$ and $L$. The unique solution of the problem {\bf (VIP)} is denoted by $x^\dagger$. 
\vskip 2mm

We need the following two lemmas to prove the convergence of Algorithm \ref{alg1}:
\vskip 2mm

\begin{lemma}{\rm \cite[Proposition 12.26]{BC2011}}\label{prox} We have
$$
\bar{x}={\rm prog}_g(z)\, \Longleftrightarrow\, \left\langle \bar{x}-z,x-\bar{x}\right\rangle \ge g(\bar{x})-g(x),\,\,\,\forall x\in \Re^m.
$$
\end{lemma}
\vskip 2mm

\begin{lemma}{\rm \cite[Corollary 2.14]{BC2011}}\label{eq} We have
$$
||\alpha x+(1-\alpha)y||^2=\alpha ||x||^2+(1-\alpha)||y||^2-\alpha(1-\alpha)||x-y||^2,\,\,\,\forall x,~y\in \Re^m,~\alpha\in \Re.
$$
\end{lemma}
\vskip 2mm

Now, we have the following first main result:
\vskip 2mm

\begin{theorem}\label{theo1}
Under the hypotheses {\rm (SP)} and {\rm (LC)}, the sequence $\left\{x_n\right\}$ generated by Algorithm \ref{alg1} converges to the unique solution $x^\dagger$ 
of the problem {\rm {\bf (VIP)}}.
\end{theorem}
\vskip 2mm

\begin{proof}
It follows from the definition of $x_{n+1}$ and Lemma \ref{prox} that 
\begin{equation}\label{h1}
\left\langle x_{n+1}-\bar{x}_n+\lambda_n F(x_n),x-x_{n+1}\right\rangle \ge \lambda_n (g(x_{n+1})-g(x)),~\forall x\in \Re^m,
\end{equation}
which, with $x=x^\dagger$, implies that 
$$
\left\langle x_{n+1}-\bar{x}_n+\lambda_n F(x_n),x^\dagger-x_{n+1}\right\rangle \ge \lambda_n (g(x_{n+1})-g(x^\dagger)).
$$
Thus we have
$$
2\left\langle x_{n+1}-\bar{x}_n,x^\dagger-x_{n+1}\right\rangle +2\lambda_n\left\langle F(x_n),x^\dagger-x_{n+1}\right\rangle \ge 2\lambda_n (g(x_{n+1})-g(x^\dagger)).
$$
This together with the equality $2\left\langle a,b\right\rangle=||a+b||^2-||a||^2+||b||^2$ and the hypothesis (SP) imply that
\begin{eqnarray}
&&||x_{n+1}-x^\dagger||^2\nonumber\\
&\le&||\bar{x}_n-x^\dagger||^2-||x_{n+1}-\bar{x}_n||^2+2\lambda_n\left\langle F(x_n),x^\dagger-x_{n+1}\right\rangle +2\lambda_n (g(x^\dagger)-g(x_{n+1}))\nonumber\\ 
&=&||\bar{x}_n-x^\dagger||^2-||x_{n+1}-\bar{x}_n||^2+2\lambda_n\left\langle F(x_n),x^\dagger-x_n\right\rangle+2\lambda_n\left\langle F(x_n),x_n-x_{n+1}\right\rangle\nonumber\\ 
&& +2\lambda_n (g(x^\dagger)-g(x_{n+1}))\nonumber\\  
&\le&||\bar{x}_n-x^\dagger||^2-||x_{n+1}-\bar{x}_n||^2-2\lambda_n\gamma||x_n-x^\dagger||^2+2\lambda_n\left\langle F(x^\dagger),x^\dagger-x_n\right\rangle\nonumber\\ 
&& +2\lambda_n\left\langle F(x_n),x_n-x_{n+1}\right\rangle+2\lambda_n (g(x^\dagger)-g(x_{n+1})).\label{h2}
\end{eqnarray}
Now, using the relation (\ref{h1}) with $n=n-1$, we obtain 
\begin{equation}\label{t1}
\left\langle x_n-\bar{x}_{n-1}+\lambda_{n-1} F(x_{n-1}),x-x_n\right\rangle \ge \lambda_{n-1} (g(x_n)-g(x)),\,\,\,\forall x\in \Re^m.
\end{equation}
Substituting $x=x_{n+1}$ into the inequality (\ref{t1}), we get 
\begin{equation}\label{h3}
\left\langle x_n-\bar{x}_{n-1}+\lambda_{n-1} F(x_{n-1}),x_{n+1}-x_n\right\rangle \ge \lambda_{n-1} (g(x_n)-g(x_{n+1})).
\end{equation}
Multiplying both sides of (\ref{h3}) by $\frac{2\lambda_n}{\lambda_{n-1}}>0$ and noting that $x_n-\bar{x}_{n-1}=\varphi(x_n-\bar{x}_n)$, we get 
\begin{equation}\label{h4}
2\frac{\varphi\lambda_n}{\lambda_{n-1}}\left\langle x_n-\bar{x}_n,x_{n+1}-x_n\right\rangle +2\lambda_n\left\langle F(x_{n-1}),x_{n+1}-x_n\right\rangle \ge 2\lambda_n (g(x_n)-g(x_{n+1})).
\end{equation}
Thus, using the identity $2\left\langle a,b\right\rangle=||a+b||^2-||a||^2+||b||^2$, we have
 the following estimate:
\begin{eqnarray}
0&\le&\frac{\varphi\lambda_n}{\lambda_{n-1}}\left[ ||x_{n+1}-\bar{x}_n||^2-||x_n-\bar{x}_n||^2-||x_{n+1}-x_n||^2\right]\nonumber\\
&&+2\lambda_n\left\langle F(x_{n-1}),x_{n+1}-x_n\right\rangle+ 2\lambda_n (g(x_{n+1})-g(x_n)).\label{h5}
\end{eqnarray}
Adding both sides of the relations (\ref{h2}) and (\ref{h5}), we obtain
\begin{eqnarray}
&&||x_{n+1}-x^\dagger||^2\nonumber\\
&\le&||\bar{x}_n-x^\dagger||^2-2\gamma\lambda_n||x_n-x^\dagger||^2-\left(1-\frac{\varphi\lambda_n}{\lambda_{n-1}}\right)||x_{n+1}-\bar{x}_n||^2\nonumber\\
&&-\frac{\varphi\lambda_n}{\lambda_{n-1}}\left[ ||x_n-\bar{x}_n||^2+||x_{n+1}-x_n||^2\right]+2\lambda_n \left\langle F(x_n)-F(x_{n-1}),x_n-x_{n+1}\right\rangle\nonumber\\
&&-2\lambda_n \left[\left\langle F(x^\dagger),x_n-x^\dagger\right\rangle +g(x_n)-g(x^\dagger)\right].\label{h6}
\end{eqnarray}
Using the Lipschitz continuity of $F$, we derive
\begin{eqnarray}
2\left\langle F(x_n)-F(x_{n-1}),x_n-x_{n+1}\right\rangle&\le &2 L||x_n-x_{n-1}|| ||x_n-x_{n+1}||\nonumber\\ 
&\le &L||x_n-x_{n-1}||^2+L ||x_n-x_{n+1}||^2.\label{h6a} 
\end{eqnarray}
Since $x^\dagger$ is a solution of the problem {\bf (VIP)}, we have
\begin{equation}\label{h6b}
\left\langle F(x^\dagger),x_n-x^\dagger\right\rangle +g(x_n)-g(x^\dagger)\ge 0.
\end{equation}
Combining the relations (\ref{h6})-(\ref{h6b}), we get
\begin{eqnarray}
&&||x_{n+1}-x^\dagger||^2\nonumber\\
&\le&||\bar{x}_n-x^\dagger||^2-2\gamma\lambda_n||x_n-x^\dagger||^2-\left(1-\frac{\varphi\lambda_n}{\lambda_{n-1}}\right)||x_{n+1}-\bar{x}_n||^2\nonumber\\
&&-\frac{\varphi\lambda_n}{\lambda_{n-1}} ||x_n-\bar{x}_n||^2+\lambda_n L||x_n-x_{n-1}||^2-\left(\frac{\varphi\lambda_n}{\lambda_{n-1}}-\lambda_n L\right)||x_{n+1}-x_n||^2.\nonumber\\
\label{h6c}
\end{eqnarray}
Moreover, from the definition of $\bar{x}_n$ and Lemma \ref{eq}, it follows that
\begin{eqnarray}\label{h7}
||x_{n+1}-x^\dagger||^2&=&\frac{\varphi}{\varphi-1}||\bar{x}_{n+1}-x^\dagger||^2-\frac{1}{\varphi-1}||\bar{x}_n-x^\dagger||^2+\frac{\varphi}{(\varphi-1)^2}||\bar{x}_{n+1}-\bar{x}_n||^2\nonumber\\
&=&\frac{\varphi}{\varphi-1}||\bar{x}_{n+1}-x^\dagger||^2-\frac{1}{\varphi-1}||\bar{x}_n-x^\dagger||^2+\frac{1}{\varphi}||x_{n+1}-\bar{x}_n||^2.
\end{eqnarray}
Combining the relations (\ref{h6c}) and (\ref{h7}), we obtain 
\begin{eqnarray}
&&\frac{\varphi}{\varphi-1}||\bar{x}_{n+1}-x^\dagger||^2\nonumber\\
&\le&
\frac{\varphi}{\varphi-1}||\bar{x}_n-x^\dagger||^2-2\gamma\lambda_n||x_n-x^\dagger||^2-\left(1+\frac{1}{\varphi}-\frac{\varphi\lambda_n}{\lambda_{n-1}}\right)||x_{n+1}-\bar{x}_n||^2\nonumber\\
&&-\frac{\varphi\lambda_n}{\lambda_{n-1}} ||x_n-\bar{x}_n||^2+\lambda_n L||x_n-x_{n-1}||^2-\left(\frac{\varphi\lambda_n}{\lambda_{n-1}}-\lambda_n L\right)||x_{n+1}-x_n||^2.\nonumber\\
\label{h8}
\end{eqnarray}
Since the sequence $\left\{\lambda_n\right\}$ is non-increasing, we obtain
\begin{equation}\label{h9}
1+\frac{1}{\varphi}-\frac{\varphi\lambda_n}{\lambda_{n-1}}\ge 1+\frac{1}{\varphi}-\varphi=0.
\end{equation}
Since $\lambda_n\to 0$, there exists $n_0\ge 1$ such that 
$$
\frac{\varphi}{\lambda_{n-1}}-2L>0,\,\,\,\forall n\ge n_0.
$$
Thus, from $\lambda_{n+1}\le \lambda_n$, it follows that
$$
\frac{\varphi\lambda_n}{\lambda_{n-1}}-\lambda_n L-\lambda_{n+1} L\ge \frac{\varphi\lambda_n}{\lambda_{n-1}}-\lambda_n L-\lambda_n L=\lambda_n \left[\frac{\varphi}{\lambda_{n-1}}-2L\right]>0,\,\,\,\forall n\ge n_0,
$$
which follows that
\begin{equation}\label{h12}
\frac{\varphi\lambda_n}{\lambda_{n-1}}-\lambda_n L>\lambda_{n+1} L,\,\,\, \forall n\ge n_0.
\end{equation}
From the relations (\ref{h8}), (\ref{h9}) and (\ref{h12}), it follows that, for all $n\ge n_0$,
\begin{eqnarray*}
\frac{\varphi}{\varphi-1}||\bar{x}_{n+1}-x^\dagger||^2&\le& \frac{\varphi}{\varphi-1}||\bar{x}_n-x^\dagger||^2-2\gamma\lambda_n||x_n-x^\dagger||^2-\frac{\varphi\lambda_n}{\lambda_{n-1}} ||x_n-\bar{x}_n||^2\nonumber\\
&&+\lambda_n L||x_n-x_{n-1}||^2-\lambda_{n+1} L||x_{n+1}-x_n||^2
\end{eqnarray*}
or
\begin{eqnarray}\label{h14}
a_{n+1}\le a_n-b_n,\,\,\,\forall n\ge n_0,
\end{eqnarray}
where 
$$
a_n:=\frac{\varphi}{\varphi-1}||\bar{x}_n-x^\dagger||^2+\lambda_n L||x_n-x_{n-1}||^2,
$$
$$
b_n=\frac{\varphi\lambda_n}{\lambda_{n-1}} ||x_n-\bar{x}_n||^2+2\gamma\lambda_n||x_n-x^\dagger||^2.
$$
Thus the sequence $\left\{a_n\right\}_{n\ge n_0}$ is non-increasing. It is obvious that $\left\{a_n\right\}_{n\ge n_0}$ is bounded from below by $0$. Thus the limit $\lim\limits_{n\to\infty}a_n$ exists and $\lim\limits_{n\to\infty}a_n\in \Re$. This implies that $\left\{\bar{x}_n\right\}$ is bounded. Thus, from the definition of $\bar{x}_n$, 
we also see that the sequence $\left\{x_n\right\}$ is bounded. This together with the fact $\lambda_n\to 0$ implies that $\lambda_n L||x_n-x_{n-1}||^2\to 0$ as 
$n\to\infty$. Hence, from the definition of $a_n$, it follows that
\begin{equation}\label{h15*}
\lim\limits_{n\to\infty}||\bar{x}_n-x^\dagger||^2\in \Re.
\end{equation}
Also, from the relation (\ref{h14}), we obtain $\sum\limits_{n=n_0}^\infty b_n <+\infty$. Thus, from the definition of $b_n$, we obtain 
$$
\mbox{\rm (S1)}\,\, \sum\limits_{n=n_0}^\infty \frac{\lambda_n}{\lambda_{n-1}} ||x_n-\bar{x}_n||^2<+\infty ,\qquad \mbox{\rm (S2)}\,\,\sum\limits_{n=n_0}^\infty \lambda_n||x_n-x^\dagger||^2<+\infty.
$$
It follows from the conditions (H3) and (S1) that $||\bar{x}_n-x_n||^2\to 0$ as $n\to\infty$. Thus, from the relation (\ref{h15*}), we get 
\begin{equation}\label{h15}
\lim\limits_{n\to\infty}||x_n-x^\dagger||^2\in \Re.
\end{equation}
On the other hand, from (H2) and (S2), we obtain that $\liminf\limits_n ||x_n-x^\dagger||^2=0$. This together with the relation (\ref{h15}) implies that 
$\lim\limits_{n\to\infty}||x_n-x^\dagger||^2=0$ or the sequence $\left\{x_n\right\}$ converges to $x^\dagger$, which solves uniquely the problem {\bf (VIP)}. This completes the proof.
\end{proof}
\vskip 2mm

\begin{remark}
It follows from (S2) that $\lambda_n||x_n-x^\dagger||^2<\frac{1}{n}$. Thus, if we choose $\lambda_n=\frac{1}{n^p}$ (with $0<p<1$), we obtain 
an estimate of the convergence rate of the sequence $\left\{x_n\right\}$ generated by Algorithm \ref{alg1} that 
$$||x_n-x^\dagger||^2=O\Big(\frac{1}{n^{1-p}}\Big).$$
\end{remark}

\section{Golden Ratio Algorithm without diminishing stepsizes}\label{main1}

In this section, we introduce a simple stepsize rule where the stepsizes will be updated over each iteration and only uses the information on the data, 
the previous approximations without the prior knowledge of Lipschitz constant. Unlike the previous section, the stepsizes generated by the next algorithm 
are bounded from below by a positive constant. Another stepsize rule for a golden ratio algorithm can be found in \cite{M2018}. For the 
sake of simplicity, we adopt the convention $\frac{0}{0}=+\infty$.
\vskip 2mm

 Now, we describe the algorithm in details as follows:
 \vskip 2mm
 
\begin{algorithm}[Golden Ratio Algorithm without diminishing stepsizes].\label{alg2}\\
\textbf{Initialization:} Choose $\bar{x}_0,~x_1,~x_0\in \Re^m$, $\lambda_0>0$, $\mu\in \left(0,\frac{\varphi}{2}\right)$\\
\textbf{Iterative Steps:} Assume that $\bar{x}_{n-1},~x_{n-1}, ~x_n$ are known, calculate $x_{n+1}$ as follows:
$$
\left \{
\begin{array}{ll}
\bar{x}_n=\frac{(\varphi-1)x_n+\bar{x}_{n-1}}{\varphi},\\
x_{n+1}=\mbox{\rm prox}_{\lambda_n g}(\bar{x}_n-\lambda_n F(x_n)),
\end{array}
\right.
$$
where
$$
\lambda_n=\min\left\{\lambda_{n-1},\frac{\mu||x_n-x_{n-1}||}{||F(x_n)-F(x_{n-1})||}\right\}.
$$
\end{algorithm}
\vskip 2mm

\begin{remark}
It follows from the definition of $\left\{\lambda_n\right\}$ that this sequence is non-increasing. Moreover, from the Lipschitz continuity of $F$ and 
in the case $F(x_n)\ne F(x_{n-1})$, we see that
$$
\frac{\mu||x_n-x_{n-1}||}{||F(x_n)-F(x_{n-1})||}\ge \frac{\mu||x_n-x_{n-1}||}{L||x_n-x_{n-1}||}=\frac{\mu}{L}.
$$
Thus we obtain by the induction that $\lambda_n\ge \min\left\{\lambda_0,\frac{\mu}{L}\right\}$ and so there exists $\lambda>0$ such that 
\begin{equation}\label{h16}
\lim_{n\to\infty}\lambda_n=\lambda>0.
\end{equation}
\end{remark}

\subsection{The convergence of Algorithm \ref{alg2}}

In this subsection, we study the convergence of Algorithm \ref{alg2}. We weaken the assumptions imposed on the cost operator $F$ where it only need to satisfy 
the above condition (LC) and the following pseudomonotone condition (PC): 
\vskip 1mm

(PC)\, $F$ is pseudomonotone, i.e.,  the following implication holds:
$$
\left\langle F(x),x-y\right\rangle \ge 0\, \Longrightarrow \,  \left\langle F(y),x-y\right\rangle \ge 0,\,\,\,\forall x,~y\in \Re^m.
$$
\vskip 2mm

We have the following second result:
\vskip 2mm

\begin{theorem}\label{theo2}
Under the conditions {\rm (LC)} and {\rm (PC)}, the sequence $\left\{x_n\right\}$ generated by Algorithm \ref{alg2} converges to a solution of the problem {\rm {\bf (VIP)}}.
\end{theorem}
\vskip 2mm

\begin{proof}
By arguing similarly to the relation (\ref{h6}) with $F$ being pseudomonotone and $x^*\in VI(F,g)$, we obtain 
\begin{eqnarray}
&&||x_{n+1}-x^*||^2\nonumber\\
&\le&||\bar{x}_n-x^*||^2-\left(1-\frac{\varphi\lambda_n}{\lambda_{n-1}}\right)||x_{n+1}-\bar{x}_n||^2-\frac{\varphi\lambda_n}{\lambda_{n-1}}\left[ ||x_n-\bar{x}_n||^2+||x_{n+1}-x_n||^2\right]\nonumber\\
&&+2\lambda_n \left\langle F(x_n)-F(x_{n-1}),x_n-x_{n+1}\right\rangle,\label{h17}
\end{eqnarray}
where with noting in (\ref{h6}) that $\left\langle F(x^*),x_n-x^*\right\rangle +g(x_n)-g(x^*)\ge 0$.
From the definition of $\lambda_n$, we see that 
\begin{eqnarray}
2\lambda_n \left\langle F(x_n)-F(x_{n-1}),x_n-x_{n+1}\right\rangle&\le&2\lambda_n ||F(x_n)-F(x_{n-1})||||x_n-x_{n+1}||\nonumber\\ 
&\le & 2\mu ||x_n-x_{n-1}||||x_n-x_{n+1}||\nonumber\\ 
&\le & \mu ||x_n-x_{n-1}||^2+\mu ||x_n-x_{n+1}||^2.\label{h19}
\end{eqnarray}
We have the following fact (see, the relation (\ref{h7})):
\begin{eqnarray}\label{h20}
||x_{n+1}-x^\dagger||^2&=&\frac{\varphi}{\varphi-1}||\bar{x}_{n+1}-x^\dagger||^2-\frac{1}{\varphi-1}||\bar{x}_n-x^\dagger||^2+\frac{1}{\varphi}||x_{n+1}-\bar{x}_n||^2.
\end{eqnarray}
Combining the relations (\ref{h17}) - (\ref{h20}), we get 
\begin{eqnarray}
&&\frac{\varphi}{\varphi-1}||\bar{x}_{n+1}-x^*||^2\nonumber\\
&\le&\frac{\varphi}{\varphi-1}||\bar{x}_n-x^*||^2-\left(1+\frac{1}{\varphi}-\frac{\varphi\lambda_n}{\lambda_{n-1}}\right)||x_{n+1}-\bar{x}_n||^2-\frac{\varphi\lambda_n}{\lambda_{n-1}} ||x_n-\bar{x}_n||^2\nonumber\\
&&+\mu ||x_n-x_{n-1}||^2-\left(\frac{\varphi\lambda_n}{\lambda_{n-1}}-\mu\right)||x_{n+1}-x_n||^2.\label{h21}
\end{eqnarray}
Since $\left\{\lambda_n\right\}$ is non-increasing, we get 
\begin{equation}\label{h22}
1+\frac{1}{\varphi}-\frac{\varphi\lambda_n}{\lambda_{n-1}}\ge 1+\frac{1}{\varphi}-\varphi=0.
\end{equation}
Since $\mu\in \left(0,\frac{\varphi}{2}\right)$, it follows from the relation (\ref{h16}) that
$$
\lim_{n\to\infty}\frac{\varphi\lambda_n}{\lambda_{n-1}}-\mu=\varphi-\mu>\mu.
$$
Therefore, there exists $n_0\ge 1$ such that
\begin{equation}\label{h23}
\frac{\varphi\lambda_n}{\lambda_{n-1}}-\mu>\mu, \,\,\,\forall n\ge n_0.
\end{equation}
From the relations (\ref{h21})-(\ref{h23}), we derive
\begin{eqnarray*}
\frac{\varphi}{\varphi-1}||\bar{x}_{n+1}&-&x^*||^2+\mu ||x_{n+1}-x_n||^2\\
&\le&\frac{\varphi}{\varphi-1}||\bar{x}_n-x^*||^2+\mu ||x_n-x_{n-1}||^2-\frac{\varphi\lambda_n}{\lambda_{n-1}} ||x_n-\bar{x}_n||^2\nonumber\\
\end{eqnarray*}
or 
\begin{eqnarray}\label{h25}
\bar{a}_{n+1}\le \bar{a}_n-\bar{b}_n,\,\,\,\forall n\ge n_0,
\end{eqnarray}
where 
$$
\bar{a}_n:=\frac{\varphi}{\varphi-1}||\bar{x}_n-x^*||^2+\mu ||x_n-x_{n-1}||^2,
$$
$$
\bar{b}_n=\frac{\varphi\lambda_n}{\lambda_{n-1}} ||x_n-\bar{x}_n||^2.
$$
Thus the limit of $\left\{\bar{a}_n\right\}$ exists and $\sum\limits_{n=n_0}^\infty \bar{b}_n<+\infty$. Hence the sequences $\left\{\bar{x}_n\right\}$ and 
$\left\{x_n\right\}$ are bounded. Morever, we also see that
$\sum\limits_{n=n_0}^\infty \frac{\varphi\lambda_n}{\lambda_{n-1}} ||x_n-\bar{x}_n||^2<+\infty$, which, together with the relation (\ref{h16}), implies that
\begin{equation}\label{h26}
\lim_{n\to\infty} ||x_n-\bar{x}_n||^2=0.
\end{equation}
Thus, since $x_n-\bar{x}_{n-1}=\varphi(x_n-\bar{x}_n)$, we obtain 
\begin{equation}\label{h27}
\lim_{n\to\infty} ||x_n-\bar{x}_{n-1}||^2=0.
\end{equation}
It is obvious that $||x_{n+1}-\bar{x}_n||^2\to 0$ as $n\to\infty$. This together with (\ref{h26}) implies that 
\begin{equation}\label{h28}
\lim_{n\to\infty} ||x_{n+1}-x_n||^2=0.
\end{equation}
From the definition of $x_{n}$ and Lemma \ref{prox}, we see that
\begin{equation}\label{h29}
\left\langle x_n-\bar{x}_{n-1}+\lambda_{n-1} F(x_{n-1}),x-x_n\right\rangle \ge \lambda_{n-1} (g(x_n)-g(x)),\,\,\,\forall x\in \Re^m.
\end{equation}
Now, assume that $p$ is a cluster point of $\left\{x_n\right\}$, i.e., there exists a subsequence $\left\{x_k\right\}$ of $\left\{x_n\right\}$ converging to $p$. 
Passing to the limit in (\ref{h29}) when $n=k\to\infty$ and using (\ref{h16}) and (\ref{h27}), we obtain 
\begin{equation}\label{h30}
\left\langle F(p),x-p\right\rangle \ge g(p)-g(x),\,\,\,\forall x\in \Re^m,
\end{equation}
which says that $p$ is a solution of the problem {\bf (VIP)}. 

Now, in order to finish the proof, we prove that the whole sequence $\left\{x_n\right\}$ converges to $p$. Indeed, assume that $\left\{x_l\right\}$ is another subsequence of $\left\{x_n\right\}$ converges to $\bar{p}\ne p$. Note that, as mentioned above, $\bar{p}$ is also a solution 
of the problem {\bf (VIP)}. Since $\lim\limits_{n\to\infty}\bar{a}_n\in \Re$ and the relation (\ref{h28}), it follows that $\lim\limits_{n\to\infty}||\bar{x}_n-x^*||^2\in \Re$ and 
thus $\lim\limits_{n\to\infty}||x_n-x^*||^2\in \Re$ for each $x^*\in VI(F,g)$. We have the following equality:
$$
2 \left\langle x_n,p-\bar{p}\right\rangle = ||x_n-\bar{p}||^2-||x_n-p||^2+||p||^2-||\bar{p}||^2.
$$
Thus, since $\lim\limits_{n\to\infty}||x_n-p||^2\in \Re$ and $\lim\limits_{n\to\infty}||x_n-\bar{p}||^2\in \Re$, we obtain that $\lim\limits_{n\to\infty}\left\langle x_n,p-\bar{p}\right\rangle\in \Re$. Set 
\begin{equation}\label{h31}
\lim\limits_{n\to\infty}\left\langle x_n,p-\bar{p}\right\rangle=M.
\end{equation}
Now, passing to the limit in (\ref{h31}) as $n=k,~l\to \infty$, we obtain 
$$
\left\langle p,p-\bar{p}\right\rangle=\lim\limits_{k\to\infty}\left\langle x_k,p-\bar{p}\right\rangle=M=\lim\limits_{l\to\infty}\left\langle x_l,p-\bar{p}\right\rangle =\left\langle \bar{p},p-\bar{p}\right\rangle.
$$
Thus $||p-\bar{p}||^2=0$ or $\bar{p}=p$. This completes the proof.
\end{proof}

\subsection{The convergence rate of Algorithm \ref{alg2}}

This subsection deals with the convergence rate of Algorithm \ref{alg2}. In order to get the rate of convergence, we choose $\mu$ in Algorithm \ref{alg2} such that 
$0<\mu<\frac{\rho}{1+\rho}\varphi$ with $\rho\in \left(0,\frac{1}{\sqrt{5}}\right)$ and assume that the operator $F$ satisfies the aforementioned conditions 
(SP) and (LC).
\vskip 2mm

Finally, we study the convergence rate of Algorithm \ref{alg2}.
\vskip 2mm

\begin{theorem}\label{theo3}
Under the conditions {\rm (SP)} and {\rm (LC)}, the sequence $\left\{x_n\right\}$ generated by Algorithm \ref{alg2} converges at least linearly to the unique solution 
$x^\dagger$ of the problem {\rm {\bf (VIP)}}.
\end{theorem}
\vskip 2mm

\begin{proof}
By arguing similarly to the relation (\ref{h8}), we have
\begin{eqnarray}
\frac{\varphi}{\varphi-1}||\bar{x}_{n+1}-x^\dagger||^2&\le& \frac{\varphi}{\varphi-1}||\bar{x}_n-x^\dagger||^2-\left(1+\frac{1}{\varphi}-\frac{\varphi\lambda_n}{\lambda_{n-1}}\right)||x_{n+1}-\bar{x}_n||^2\nonumber\\
&&+\mu ||x_n-x_{n-1}||^2-\frac{\varphi\lambda_n}{\lambda_{n-1}} ||x_n-\bar{x}_n||^2\nonumber\\
&&-\left(\frac{\varphi\lambda_n}{\lambda_{n-1}}-\mu\right)||x_{n+1}-x_n||^2-2\gamma\lambda_n||x_n-x^\dagger||^2.\label{h32}
\end{eqnarray}
Thus we have
\begin{eqnarray}
\frac{\varphi}{\varphi-1}||\bar{x}_{n+1}&-& x^\dagger||^2+\left(\frac{\varphi\lambda_n}{\lambda_{n-1}}-\mu\right)||x_{n+1}-x_n||^2\nonumber\\
&\le&\frac{\varphi}{\varphi-1}||\bar{x}_n-x^\dagger||^2+\mu ||x_n-x_{n-1}||^2\nonumber\\
&&-\left(1+\frac{1}{\varphi}-\frac{\varphi\lambda_n}{\lambda_{n-1}}\right)||x_{n+1}-\bar{x}_n||^2-2\gamma\lambda_n||x_n-x^\dagger||^2.\label{h32}
\end{eqnarray}
Note that, from the definition of $\bar{x}_n$, we obtain
$$
x_n=\frac{\varphi}{\varphi-1}\bar{x}_n-\frac{1}{\varphi-1}\bar{x}_{n-1}.
$$
Thus it follows from Lemma \ref{eq} that
\begin{equation}\label{h33}
||x_n-x^\dagger||^2=\frac{\varphi}{\varphi-1}||\bar{x}_n-x^\dagger||^2-\frac{1}{\varphi-1}||\bar{x}_{n-1}-x^\dagger||^2+\frac{\varphi}{(\varphi-1)^2}||\bar{x}_n-\bar{x}_{n-1}||^2.
\end{equation}
From the relations (\ref{h32}) and (\ref{h33}), we see that
\begin{eqnarray}
&&\frac{\varphi}{\varphi-1}||\bar{x}_{n+1}-x^\dagger||^2+\left(\frac{\varphi\lambda_n}{\lambda_{n-1}}-\mu\right)||x_{n+1}-x_n||^2\nonumber\\
&\le &\frac{\varphi}{\varphi-1}(1-2\gamma\lambda_n)||\bar{x}_n-x^\dagger||^2
+\frac{2\gamma\lambda_n}{\varphi-1}||\bar{x}_{n-1}-x^\dagger||^2+\mu ||x_n-x_{n-1}||^2\nonumber\\&&-\left(1+\frac{1}{\varphi}-\frac{\varphi\lambda_n}{\lambda_{n-1}}\right)||x_{n+1}-\bar{x}_n||^2-\frac{2\gamma\lambda_n\varphi}{(\varphi-1)^2}||\bar{x}_n-\bar{x}_{n-1}||^2\nonumber\\
&\le& \frac{\varphi}{\varphi-1}(1-2\gamma\lambda_n)||\bar{x}_n-x^\dagger||^2+\frac{2\gamma\lambda_n}{\varphi-1}||\bar{x}_{n-1}-x^\dagger||^2+\mu ||x_n-x_{n-1}||^2, \label{h34}
\end{eqnarray}
in which the last inequality follows from the following:
$$
1+\frac{1}{\varphi}-\frac{\varphi\lambda_n}{\lambda_{n-1}}\ge 1+\frac{1}{\varphi}-\varphi \ge 0,
$$ 
$$
\frac{2\gamma\lambda_n\varphi}{(\varphi-1)^2}\ge 0.
$$ 
Now, if we set $a_n=\frac{\varphi}{\varphi-1}||\bar{x}_n-x^\dagger||^2$, then the inequality (\ref{h34}) can be rewritten as follows:
\begin{eqnarray}
a_{n+1}&+&\left(\frac{\varphi\lambda_n}{\lambda_{n-1}}-\mu\right)||x_{n+1}-x_n||^2\nonumber\\
&\le&(1-2\gamma\lambda_n)a_n+\frac{2\gamma\lambda_n}{\varphi}a_{n-1}+\mu ||x_n-x_{n-1}||^2.\label{h35}
\end{eqnarray}
Let $\beta\in (1/\varphi,1)$ be fixed. From the relation (\ref{h16}) and the fact $0<\mu<\frac{\rho\varphi}{1+\rho}$, we have
$$
\lim_{n\to \infty }\frac{\varphi\lambda_n}{\lambda_{n-1}}-\mu=\varphi-\mu>\frac{\mu}{\rho},\quad \lim_{n\to \infty }\frac{\lambda_n}{\varphi\lambda}=\frac{1}{\varphi}<\beta<1.
$$
Thus there exists $n_0\ge 1$ such that 
\begin{equation}\label{h36}
\frac{\varphi\lambda_n}{\lambda_{n-1}}-\mu>\frac{\mu}{\rho},\quad  \frac{\lambda_n}{\varphi\lambda}<\beta,\,\,\,\forall n\ge n_0.
\end{equation}
Combining the relations (\ref{h35}) and (\ref{h36}) and noting that $\lambda_n\ge \lambda$, we obtain 
\begin{eqnarray}
a_{n+1}+\frac{\mu}{\rho}||x_{n+1}-x_n||^2&\le&(1-2\gamma\lambda)a_n+2\gamma\lambda \beta a_{n-1}+\mu ||x_n-x_{n-1}||^2.\label{h37}
\end{eqnarray}
Set $b_n=\frac{\mu}{\rho}||x_{n}-x_{n-1}||^2$ and $\alpha=2\gamma\lambda$. Then we can rewrite the relation (\ref{h37}) as follows:
\begin{eqnarray}
a_{n+1}+b_{n+1}&\le&(1-\alpha)a_n+\alpha \beta a_{n-1}+\rho b_n.\label{h38}
\end{eqnarray}
On the other hand, since $F$ is strongly monotone, we get
\begin{eqnarray*}
\gamma||x_n-x_{n-1}||^2&\le&\left\langle F(x_n)-F(x_{n-1}),x_n-x_{n-1}\right\rangle \\
&\le& ||F(x_n)-F(x_{n-1})|| ||x_n-x_{n-1}||\\
&\le& \frac{\mu ||x_n-x_{n-1}||^2}{\lambda_n}. 
\end{eqnarray*}
Thus we have $\gamma\lambda_n\le \mu$ or $\alpha=2\gamma\lambda_n\le 2\mu <\frac{2\rho\varphi}{1+\rho}<1$. Let $r_1>0$ and $r_2>0$. Now, we can rewrite the 
relation (\ref{h38}) as follows:
\begin{eqnarray}
a_{n+1}+r_1 a_n+b_{n+1}&\le& r_2(a_n+r_1 a_{n-1})+\rho b_n\nonumber\\
&&+(1-\alpha-r_2+r_1)a_n+(\alpha \beta -r_1 r_2) a_{n-1}.\label{h39}
\end{eqnarray}
Choose $r_1>0$ and $r_2>0$ such that $1-\alpha-r_2+r_1=0$ and $\alpha \beta -r_1 r_2=0$. Thus we have
\begin{equation}\label{h40}
r_1=\frac{\alpha-1+\sqrt{(\alpha-1)^2+4\alpha\beta}}{2},\quad r_2=\frac{1-\alpha+\sqrt{(\alpha-1)^2+4\alpha\beta}}{2}.
\end{equation}
Consider the function 
$$
f(t)=\frac{1-t+\sqrt{(t-1)^2+4 t\beta}}{2},\,\,\, \forall t\in (0,1).
$$
Then we have 
$$
f'(t)=\frac{-1+\frac{t-1+2\beta}{\sqrt{(t-1)^2+4 t\beta}}}{2}=\frac{4\beta(\beta-1)}{2\sqrt{(t-1)^2+4 t\beta}\left(t-1+2\beta+\sqrt{(t-1)^2+4 t\beta}\right)}<0
$$
because of $\frac{1}{\varphi}<\beta<1$. Thus $f(t)$ is non-increasing on $(0,1)$. Hence we have $0<r_2=f(\alpha)<f(0)=1$. Now, set $\theta=\max \left\{\rho,r_2\right\}$ 
and note $\theta\in (0,1)$. Then it follows from the relation (\ref{h39}) that
\begin{eqnarray}
a_{n+1}+r_1 a_n+b_{n+1}&\le& \theta (a_n+r_1 a_{n-1}+ b_n),\,\,\,\forall n\ge n_0.\label{h41}
\end{eqnarray}
Thus, by induction, we obtain  
\begin{eqnarray}
a_{n+1}+r_1 a_n+b_{n+1} \le \theta^{n-n_0+1} (a_{n_0}+r_1 a_{n_0-1}+ b_{n_0}),\,\,\,~\forall n\ge n_0\label{h41}
\end{eqnarray}
and so we can reduce that 
$$
\frac{\varphi}{\varphi-1}||x_{n+1}-x^\dagger||^2=a_{n+1}\le \theta^{n-n_0+1} (a_{n_0}+r_1 a_{n_0-1}+ b_{n_0}),\,\,\,\forall n\ge n_0,
$$
or
$$
||x_{n+1}-x^\dagger||^2\le M\theta^n,\,\,\,\forall n\ge n_0,
$$
where $M=\frac{(\varphi-1)(a_{n_0}+r_1 a_{n_0-1}+ b_{n_0})}{\varphi\theta^{n_0-1}}$. This completes the proof.
\end{proof}

\section{Numerical experiments}\label{example}

In this section, we perform several experiments to show the numerical behaviour of the proposed algorithm (Agorithm \ref{alg2}) in comparison with 
other algorithms. All the programs are written in Matlab 7.0 and computed on a PC Desktop Intel(R) Core(TM) i5-3210M CPU @ 2.50GHz, RAM 2.00 GB.
\vskip 2mm

\noindent
\textbf{Example 4.1} In this example, our problem of interest is a sparse logistic regression which is a popular problem in machine learning applications: 
$$
\min_{x\in \Re^m} J(x)=\sum_{i=1}^N\log\left(1+\exp(-b_i\left\langle a_i,x\right\rangle)\right)+\gamma ||x||_1,
$$
where $x\in \Re^m$, $a_i\in \Re^m$, $b_i\in \left\{-1,1\right\}$, $\gamma>0$.
\vskip 1mm

\begin{figure}[!ht]
\begin{minipage}[b]{0.45\textwidth}
\centering
\includegraphics[height=5cm,width=6cm]{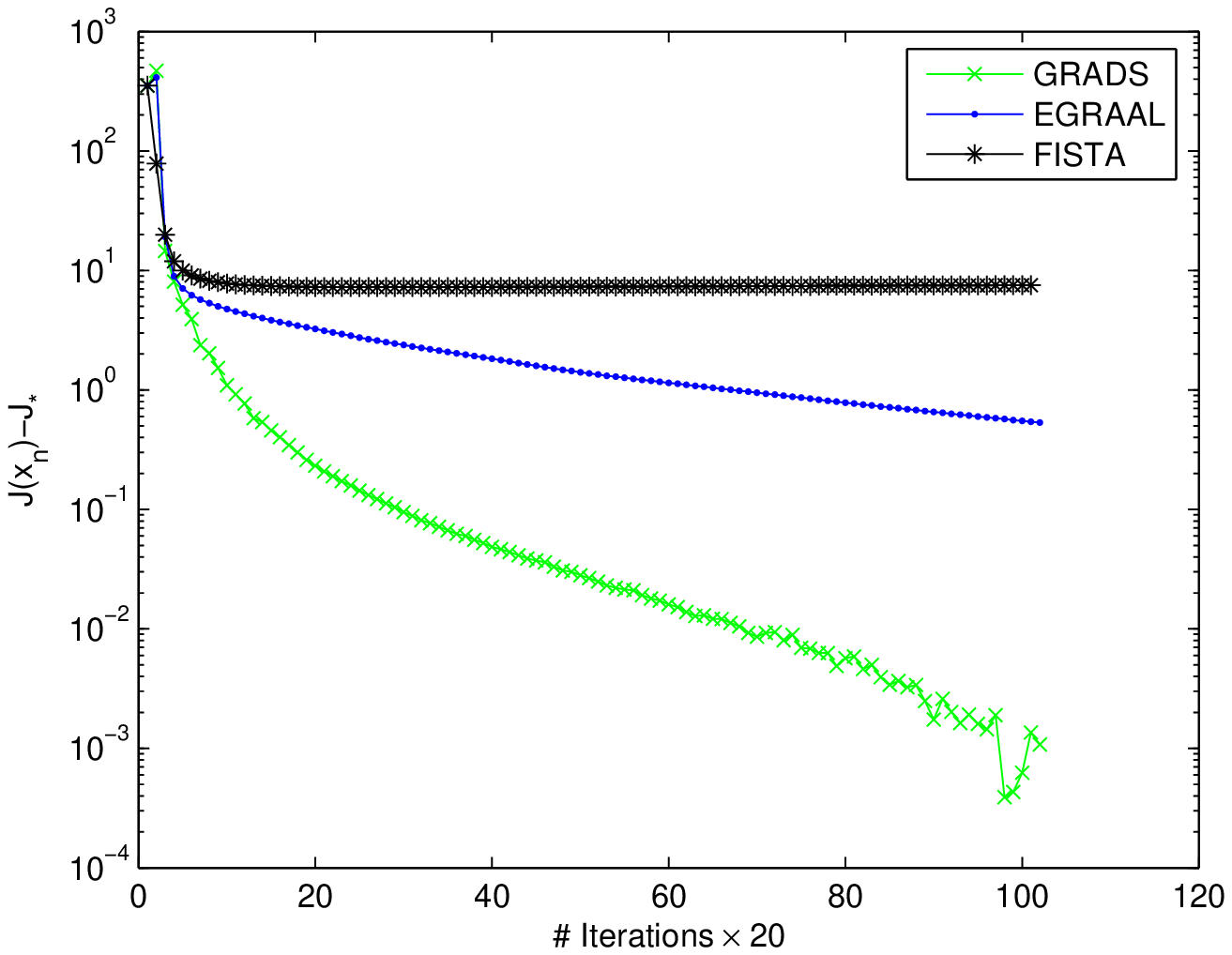}
\caption{Example 1 for $N=100$ and $m=300$}\label{fig1}
\end{minipage}
\hfill
\begin{minipage}[b]{0.45\textwidth}
\centering
\includegraphics[height=5cm,width=6cm]{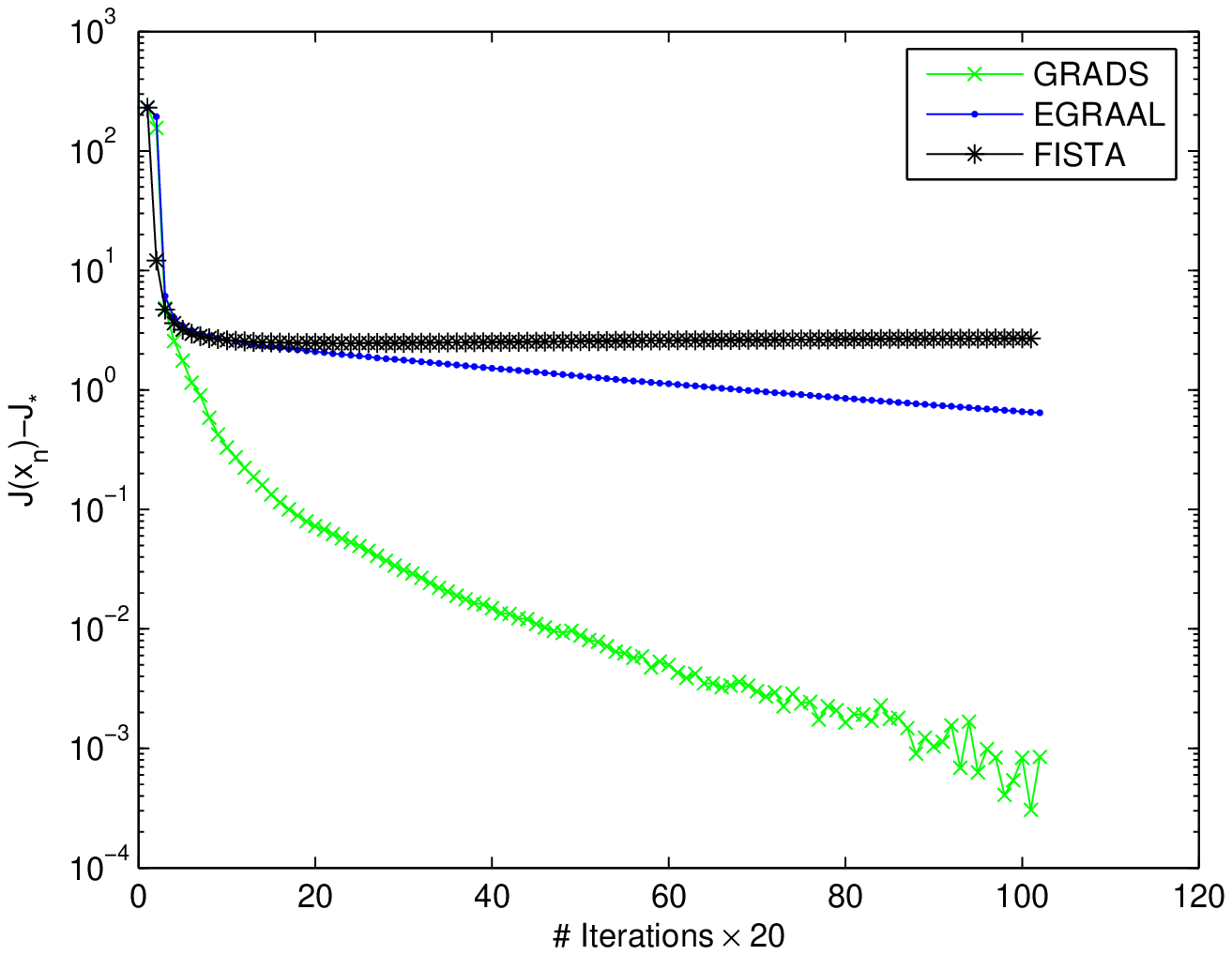}
\caption{Example 1 for $N=200$ and $m=500$}\label{fig2}
\end{minipage}
\end{figure}

 Let $K$ be a matrix of size $N\times m$ defined by $K_{ij}=-b_ia_{ij}$ 
and $\bar{f}(y)=\sum_{i=1}^N\log\left(1+\exp(y_i)\right)$. Then, the objective function in our problem is $J(x)=f(x)+g(x)$ with $f(x)=\bar{f}(Kx)$ and 
$g(x)=\gamma ||x||_1$. This problem is equivalent the considered problem with $F=\nabla f$. We compare Algorithm \ref{alg2} (GRADS) with EGRAAL 
in \cite{M2018} and FISTA with constant stepsize in \cite{BT2009}. We do not include the algorithms with linesearch procedures because they require 
many computations over each iteration which is time-consuming. Note that the algorithm FISTA requires the Lipschitz constant of $F$ ($L_{\nabla f}=||K^TK||/4$) 
while other algorithms are not. All entries of $a_i$ and $b_i$ are generated randomly and we choose $\gamma=0.005 ||A^Tb||_{\infty}$. We choose $\mu=0.45 \varphi$, 
$\lambda_0=1$ for Algorithm \ref{alg2} (GRADS); $\lambda_0=\bar{\lambda}=1$, $\phi=\varphi$ for the algorithm ERGAAL. The starting points are generated 
randomly in $(0,1]$. The mapping \textit{prox} is computed by the function \textit{fmincon} in Matlab. The results are shown on Figure \ref{fig1} and Figure \ref{fig2}. The execution times for the algorithms are almost equivalent. 
In these figures, $J_*$ is the most minimum value of $J$ generated by all the algorithms with the stopping criterion $||x_n-{\rm prox}_{g}(x_n-\nabla f(x_n))||\le 10^{-3}$.

\begin{figure}[!ht]
\begin{minipage}[b]{0.45\textwidth}
\centering
\includegraphics[height=5cm,width=6cm]{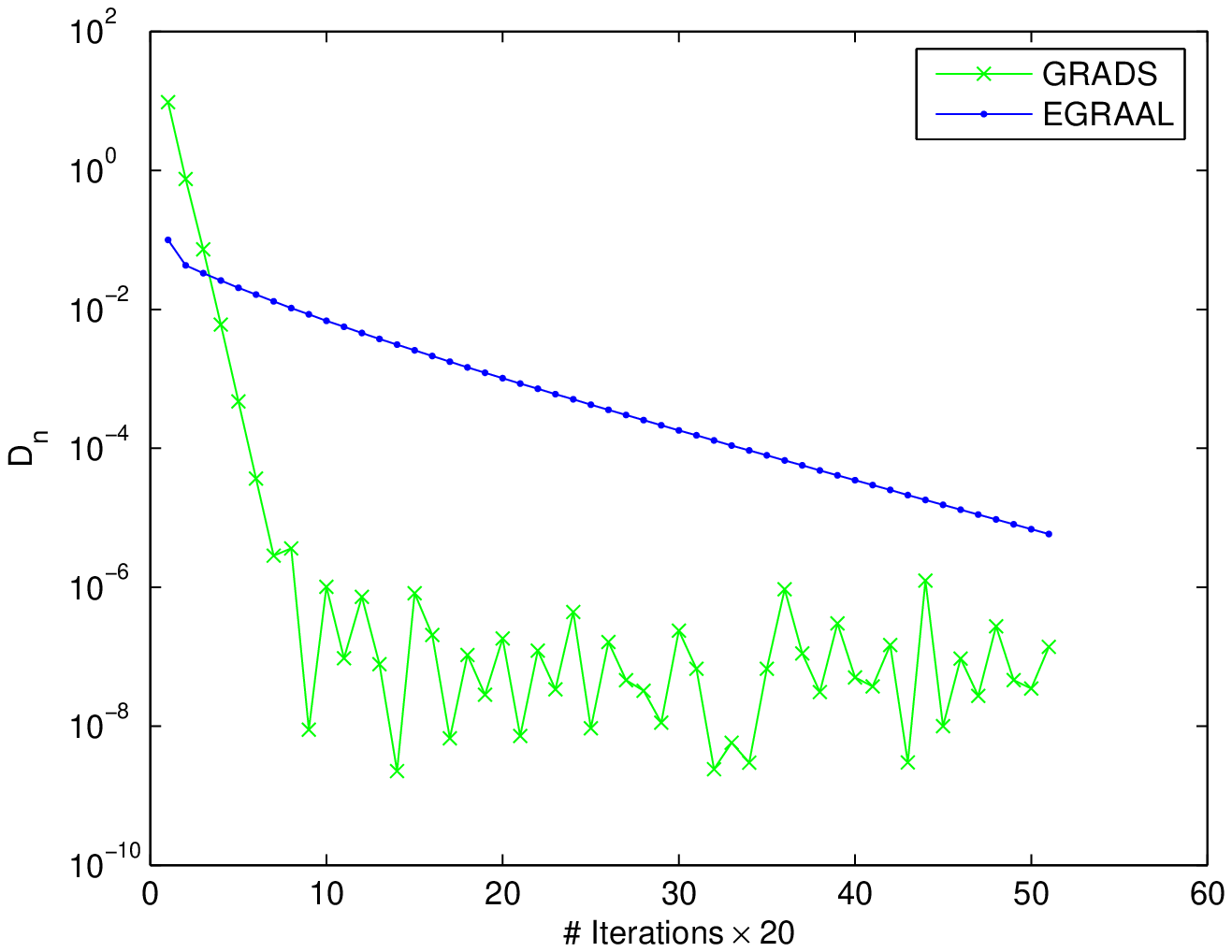}
\caption{Example 2 for $m=300$}\label{fig3}
\end{minipage}
\hfill
\begin{minipage}[b]{0.45\textwidth}
\centering
\includegraphics[height=5cm,width=6cm]{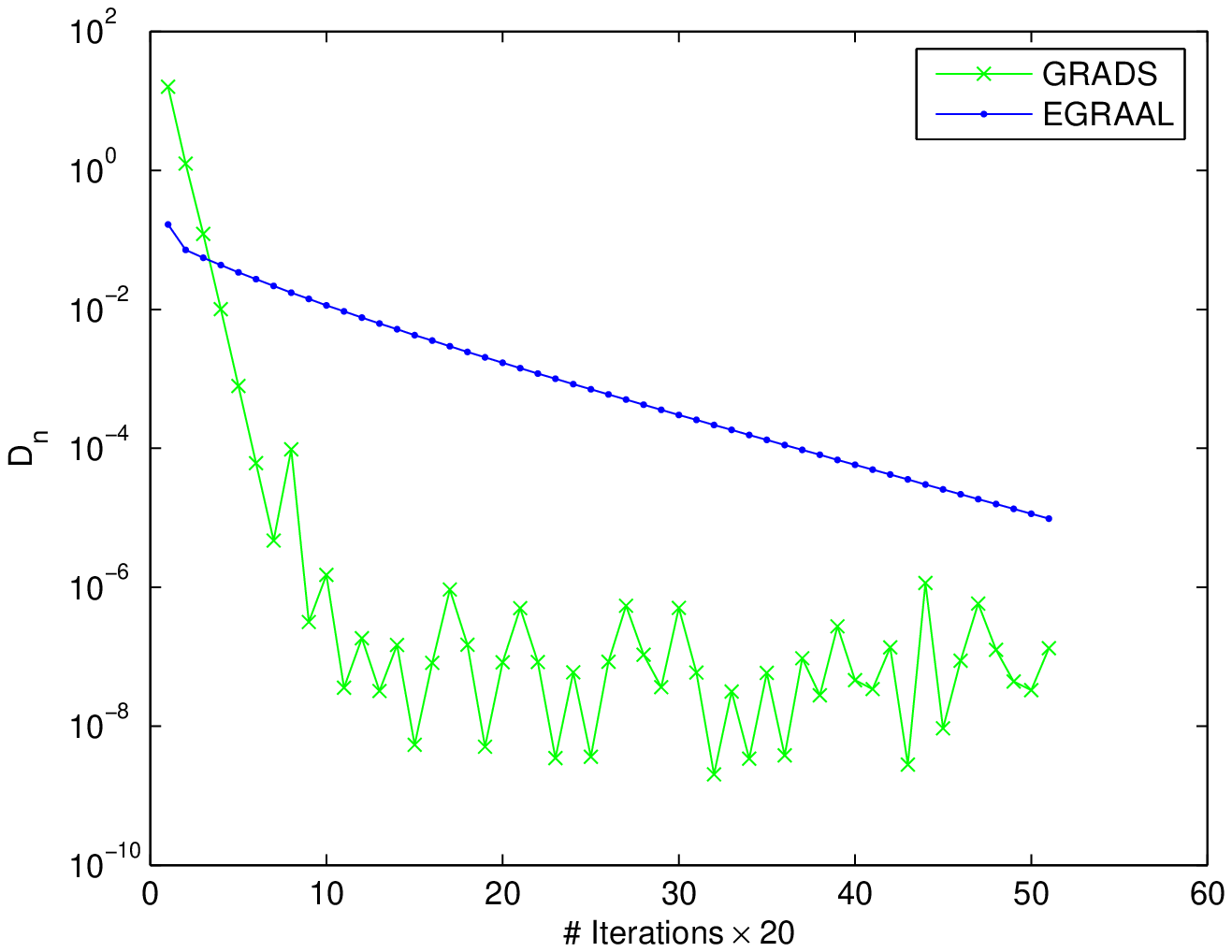}
\caption{Example 2 for $m=500$}\label{fig4}
\end{minipage}
\end{figure}

\noindent \textbf{Example 4.2} Consider the nonlinear problem, which is presented by Sun  \cite{S1994}, for the operator $A:\Re^m\to\Re^m$ of 
the form:
$$
\begin{cases}
F(x)=F_1(x)+F_2(x),\\ 
 F_1(x)=(f_1(x),f_2(x),\cdots,f_m(x)),\\
 F_2(x)=Dx+c,\\
f_i(x)=x_{i-1}^2+x_i^2+x_{i-1}x_i+x_ix_{i+1},~i=1,2,\cdots,m,\\
x_0=x_{m+1}=0,
\end{cases}
$$
where $c=(-1,-1,\cdots,-1)$ and $D$ is a square matrix of order $m$, which is given by 
$$
d_{ij}=
\begin{cases}
4,&\mbox{if}\quad i=j,\\ 
1,&\mbox{if}\quad i-j=1,\\ 
-2,&\mbox{if}\quad i-j=-1,\\ 
0,&\mbox{otherwise}.
\end{cases}
$$
The feasible set is $C=\Re^m_+$. This problem is equivalent to the considered problem with $g=\delta_C$. In this case, the mapping \textit{prox} is the projection 
on the set $C$ and it is computed by the function \textit{quadprog} in Matlab 7.0. Since the Lipschitz constant of $F$ is unknown, we do not include the comparison 
with the algorithm FISTA. We use the sequence $D_n=||x_{n+1}-\bar{x}_n||^2+||\bar{x}_n-x_n||^2$ for each $n=0,1,2,\cdots$ to compare the computational performance of 
the algorithms. Figure \ref{fig3} and Figure \ref{fig4} describe the results in this example.

The numerical results here have illustrated that the proposed algorithm works well and also has competitive advantage over other algorithms.
\section{Conclusions}
In this paper, we have introduced the two golden ratio algorithms with two simple stepsize rules for solving pseudomonotone 
and Lipschitz variational inequalities in finite dimensional Hilbert spaces. The first algorithm uses a sequence of stepsizes taken priorly with some suitable 
properties while the second one itself generates variable stepsizes which are explicitly computed in each iteration and without a linesearch procedure to be run. 
We have established the convergence as well as the convergence rate of the new algorithms under appropriate conditions. The theoretical results have been 
illustrated by some our numerical experiments.

Our results can be extended to many promising directions, such as multi-value variational inequalities \cite{FC2014}, 
equilibrium problems \cite{MO1992,BO1994,HCX18,AH2018}, problem {\bf (VIP)} incorporated with 
fixed point problems, systems of variational inequalities and mixed equilibrium problems \cite{C2008,H2017COAP,HAM2017,Y2009,Y2010}, 
weak and strong convergence in Hilbert spaces as well as extensions to Banach spaces \cite{KRS2011,HS18}. This is surely our future goals.
 
\section*{Acknowledgement}
The authors would like to thank the Associate Editor and the two anonymous referees 
for their valuable comments and suggestions which helped us very much in improving the original version of this paper. 
The research of the first author was supported by the National Foundation for Science and Technology Development (NAFOS-TED) of 
Vietnam under grant number 101.01-2017.315. The research work was also supported by the National Natural Science Foundation 
of China (11771067) and the Applied Basic Project of Sichuan Province (19YYJC0157).
We also would like to thank \textbf{Dr. Yura Malitsky} for sending us the paper \cite{M2018}.

\end{document}